\documentclass[10.5pt]{article}
\usepackage{amsmath}
\usepackage{amssymb}
\usepackage{amscd}
\usepackage{xspace}
\usepackage{verbatim}
\newcommand{\mysection}[1]{
\section{#1}\setcounter{equation}{0}}
\title{\bf Maximal solutions of nonlinear parabolic equations with absorption}
\author{
 {\bf Laurent V\'eron}\\[2mm]
{\small Laboratoire de Math\'ematiques et Physique Th\'eorique, }\\
{\small  Universit\'e Fran\c{c}ois Rabelais,  Tours,  FRANCE}}

\date{}
\begin{document}
\maketitle
{\small {\bf Abstract} We study the existence and the uniqueness  of the solution of the problem (P):  $\partial_tu-\Delta u+f(u)=0$ in $Q:=\Omega\times (0,\infty)$, $u=\infty$ on the parabolic boundary $\partial_{p}Q$ when $\Omega$ is a  domain in $\mathbb R^N$ with a compact boundary and $f$ a continuous increasing function satisfying super linear growth condition.  We prove that in most cases, the existence and uniqueness is reduced to the same property for the associated stationary equation in $\Omega$.
}

\vspace{1mm}
\hspace{.05in}

\newcommand{\txt}[1]{\;\text{ #1 }\;}
\newcommand{\tbf}{\textbf}
\newcommand{\tit}{\textit}
\newcommand{\tsc}{\textsc}
\newcommand{\trm}{\textrm}
\newcommand{\mbf}{\mathbf}
\newcommand{\mrm}{\mathrm}
\newcommand{\bsym}{\boldsymbol}
\newcommand{\scs}{\scriptstyle}
\newcommand{\sss}{\scriptscriptstyle}
\newcommand{\txts}{\textstyle}
\newcommand{\dsps}{\displaystyle}
\newcommand{\fnz}{\footnotesize}
\newcommand{\scz}{\scriptsize}
\newcommand{\be}{
\begin{equation}
}
\newcommand{\bel}[1]{
\begin{equation}
\label{#1}}
\newcommand{\ee}{
\end{equation}
}
\newcommand{\eqnl}[2]{
\begin{equation}
\label{#1}{#2}
\end{equation}
}
\newtheorem{subn}{\name}
\renewcommand{\thesubn}{}
\newcommand{\bsn}[1]{\def\name{#1}
\begin{subn}}
\newcommand{\esn}{
\end{subn}}
\newtheorem{sub}{\name}[section]
\newcommand{\dn}[1]{\def\name{#1}}   
\newcommand{\bs}{
\begin{sub}}
\newcommand{\es}{
\end{sub}}
\newcommand{\bsl}[1]{
\begin{sub}\label{#1}}
\newcommand{\bth}[1]{\def\name{Theorem}
\begin{sub}\label{t:#1}}
\newcommand{\blemma}[1]{\def\name{Lemma}
\begin{sub}\label{l:#1}}
\newcommand{\bcor}[1]{\def\name{Corollary}
\begin{sub}\label{c:#1}}
\newcommand{\bdef}[1]{\def\name{Definition}
\begin{sub}\label{d:#1}}
\newcommand{\bprop}[1]{\def\name{Proposition}
\begin{sub}\label{p:#1}}
\newcommand{\R}{\eqref}
\newcommand{\rth}[1]{Theorem~\ref{t:#1}}
\newcommand{\rlemma}[1]{Lemma~\ref{l:#1}}
\newcommand{\rcor}[1]{Corollary~\ref{c:#1}}
\newcommand{\rdef}[1]{Definition~\ref{d:#1}}
\newcommand{\rprop}[1]{Proposition~\ref{p:#1}}
\newcommand{\BA}{
\begin{array}}
\newcommand{\EA}{
\end{array}}
\newcommand{\BAN}{\renewcommand{\arraystretch}{1.2}
\setlength{\arraycolsep}{2pt}
\begin{array}}
\newcommand{\BAV}[2]{\renewcommand{\arraystretch}{#1}
\setlength{\arraycolsep}{#2}
\begin{array}}
\newcommand{\BSA}{
\begin{subarray}}
\newcommand{\ESA}{
\end{subarray}}
\newcommand{\BAL}{
\begin{aligned}}
\newcommand{\EAL}{
\end{aligned}}
\newcommand{\BALG}{
\begin{alignat}}
\newcommand{\EALG}{
\end{alignat}}
\newcommand{\BALGN}{
\begin{alignat*}}
\newcommand{\EALGN}{
\end{alignat*}}
\newcommand{\note}[1]{\textit{#1.}\hspace{2mm}}
\newcommand{\Proof}{\note{Proof}}
\newcommand{\qeda}{\hspace{10mm}\hfill $\square$}
\newcommand{\qed}{\\
${}$ \hfill $\square$}
\newcommand{\Remark}{\note{Remark}}
\newcommand{\modin}{$\,$\\
[-4mm] \indent}
\newcommand{\forevery}{\quad \forall}
\newcommand{\set}[1]{\{#1\}}
\newcommand{\setdef}[2]{\{\,#1:\,#2\,\}}
\newcommand{\setm}[2]{\{\,#1\mid #2\,\}}
\newcommand{\lra}{\longrightarrow}
\newcommand{\lla}{\longleftarrow}
\newcommand{\llra}{\longleftrightarrow}
\newcommand{\Lra}{\Longrightarrow}
\newcommand{\Lla}{\Longleftarrow}
\newcommand{\Llra}{\Longleftrightarrow}
\newcommand{\warrow}{\rightharpoonup}
\newcommand{
\paran}[1]{\left (#1 \right )}
\newcommand{\sqbr}[1]{\left [#1 \right ]}
\newcommand{\curlybr}[1]{\left \{#1 \right \}}
\newcommand{\abs}[1]{\left |#1\right |}
\newcommand{\norm}[1]{\left \|#1\right \|}
\newcommand{
\paranb}[1]{\big (#1 \big )}
\newcommand{\lsqbrb}[1]{\big [#1 \big ]}
\newcommand{\lcurlybrb}[1]{\big \{#1 \big \}}
\newcommand{\absb}[1]{\big |#1\big |}
\newcommand{\normb}[1]{\big \|#1\big \|}
\newcommand{
\paranB}[1]{\Big (#1 \Big )}
\newcommand{\absB}[1]{\Big |#1\Big |}
\newcommand{\normB}[1]{\Big \|#1\Big \|}

\newcommand{\thkl}{\rule[-.5mm]{.3mm}{3mm}}
\newcommand{\thknorm}[1]{\thkl #1 \thkl\,}
\newcommand{\trinorm}[1]{|\!|\!| #1 |\!|\!|\,}
\newcommand{\bang}[1]{\langle #1 \rangle}
\def\angb<#1>{\langle #1 \rangle}
\newcommand{\vstrut}[1]{\rule{0mm}{#1}}
\newcommand{\rec}[1]{\frac{1}{#1}}
\newcommand{\opname}[1]{\mbox{\rm #1}\,}
\newcommand{\supp}{\opname{supp}}
\newcommand{\dist}{\opname{dist}}
\newcommand{\myfrac}[2]{{\displaystyle \frac{#1}{#2} }}
\newcommand{\myint}[2]{{\displaystyle \int_{#1}^{#2}}}
\newcommand{\mysum}[2]{{\displaystyle \sum_{#1}^{#2}}}
\newcommand {\dint}{{\displaystyle \int\!\!\int}}
\newcommand{\q}{\quad}
\newcommand{\qq}{\qquad}
\newcommand{\hsp}[1]{\hspace{#1mm}}
\newcommand{\vsp}[1]{\vspace{#1mm}}
\newcommand{\ity}{\infty}
\newcommand{\prt}{\partial}
\newcommand{\sms}{\setminus}
\newcommand{\ems}{\emptyset}
\newcommand{\ti}{\times}
\newcommand{\pr}{^\prime}
\newcommand{\ppr}{^{\prime\prime}}
\newcommand{\tl}{\tilde}
\newcommand{\sbs}{\subset}
\newcommand{\sbeq}{\subseteq}
\newcommand{\nind}{\noindent}
\newcommand{\ind}{\indent}
\newcommand{\ovl}{\overline}
\newcommand{\unl}{\underline}
\newcommand{\nin}{\not\in}
\newcommand{\pfrac}[2]{\genfrac{(}{)}{}{}{#1}{#2}}

\def\ga{\alpha}     \def\gb{\beta}       \def\gg{\gamma}
\def\gc{\chi}       \def\gd{\delta}      \def\ge{\epsilon}
\def\gth{\theta}                         \def\vge{\varepsilon}
\def\gf{\phi}       \def\vgf{\varphi}    \def\gh{\eta}
\def\gi{\iota}      \def\gk{\kappa}      \def\gl{\lambda}
\def\gm{\mu}        \def\gn{\nu}         \def\gp{\pi}
\def\vgp{\varpi}    \def\gr{\rho}        \def\vgr{\varrho}
\def\gs{\sigma}     \def\vgs{\varsigma}  \def\gt{\tau}
\def\gu{\upsilon}   \def\gv{\vartheta}   \def\gw{\omega}
\def\gx{\xi}        \def\gy{\psi}        \def\gz{\zeta}
\def\Gg{\Gamma}     \def\Gd{\Delta}      \def\Gf{\Phi}
\def\Gth{\Theta}
\def\Gl{\Lambda}    \def\Gs{\Sigma}      \def\Gp{\Pi}
\def\Gw{\Omega}     \def\Gx{\Xi}         \def\Gy{\Psi}

\def\CS{{\mathcal S}}   \def\CM{{\mathcal M}}   \def\CN{{\mathcal N}}
\def\CR{{\mathcal R}}   \def\CO{{\mathcal O}}   \def\CP{{\mathcal P}}
\def\CA{{\mathcal A}}   \def\CB{{\mathcal B}}   \def\CC{{\mathcal C}}
\def\CD{{\mathcal D}}   \def\CE{{\mathcal E}}   \def\CF{{\mathcal F}}
\def\CG{{\mathcal G}}   \def\CH{{\mathcal H}}   \def\CI{{\mathcal I}}
\def\CJ{{\mathcal J}}   \def\CK{{\mathcal K}}   \def\CL{{\mathcal L}}
\def\CT{{\mathcal T}}   \def\CU{{\mathcal U}}   \def\CV{{\mathcal V}}
\def\CZ{{\mathcal Z}}   \def\CX{{\mathcal X}}   \def\CY{{\mathcal Y}}
\def\CW{{\mathcal W}} \def\CQ{{\mathcal Q}}
\def\BBA {\mathbb A}   \def\BBb {\mathbb B}    \def\BBC {\mathbb C}
\def\BBD {\mathbb D}   \def\BBE {\mathbb E}    \def\BBF {\mathbb F}
\def\BBG {\mathbb G}   \def\BBH {\mathbb H}    \def\BBI {\mathbb I}
\def\BBJ {\mathbb J}   \def\BBK {\mathbb K}    \def\BBL {\mathbb L}
\def\BBM {\mathbb M}   \def\BBN {\mathbb N}    \def\BBO {\mathbb O}
\def\BBP {\mathbb P}   \def\BBR {\mathbb R}    \def\BBS {\mathbb S}
\def\BBT {\mathbb T}   \def\BBU {\mathbb U}    \def\BBV {\mathbb V}
\def\BBW {\mathbb W}   \def\BBX {\mathbb X}    \def\BBY {\mathbb Y}
\def\BBZ {\mathbb Z}

\def\GTA {\mathfrak A}   \def\GTB {\mathfrak B}    \def\GTC {\mathfrak C}
\def\GTD {\mathfrak D}   \def\GTE {\mathfrak E}    \def\GTF {\mathfrak F}
\def\GTG {\mathfrak G}   \def\GTH {\mathfrak H}    \def\GTI {\mathfrak I}
\def\GTJ {\mathfrak J}   \def\GTK {\mathfrak K}    \def\GTL {\mathfrak L}
\def\GTM {\mathfrak M}   \def\GTN {\mathfrak N}    \def\GTO {\mathfrak O}
\def\GTP {\mathfrak P}   \def\GTR {\mathfrak R}    \def\GTS {\mathfrak S}
\def\GTT {\mathfrak T}   \def\GTU {\mathfrak U}    \def\GTV {\mathfrak V}
\def\GTW {\mathfrak W}   \def\GTX {\mathfrak X}    \def\GTY {\mathfrak Y}
\def\GTZ {\mathfrak Z}   \def\GTQ {\mathfrak Q}

\font\Sym= msam10 
\def\SYM#1{\hbox{\Sym #1}}
\newcommand{\bdw}{\prt\Gw\xspace}
\medskip
\mysection {Introduction}
Let  $\Gw$ be a  bounded domain in $\BBR^N$ with  boundary $\prt\Gw:=\Gg$, $Q^\Gw_{T}:=\Gw\ti (0,T)$ ($0<T\leq \infty$) and $\prt_{p}Q=\overline\Gw\ti{0}\cup \prt\Gw\ti (0,T]$. We denote by $\gr_{_{\partial\Gw}}(x)$ the distance from $x$ to $\prt\Gw$ and by 
$d_{_{P}}(x,t)=\min\{\gr_{_{\partial\Gw}}(x),t\}$ the product distance from $(x,t)\in Q^\Gw_{\infty}$ to $\prt_{p}Q^\Gw_{\infty}$ . If 
$f\in C(\BBR)$, we say that a function $u\in C^{2,1}(Q^\Gw_{\infty})$  solution of 
\begin{equation}\label{A1}
u_{t}-\Gd u+f(u)=0,
\end{equation}
in $Q^\Gw_{\infty}$ is a large solution of (\ref{A1}) in $Q^\Gw_{\infty}$ if it satisfies
\begin{equation}\label{A2}
\lim_{d_{_{P}}(x,t)\to 0}u(x,t)=\infty.
\end{equation}
The existence of such a $u$ is associated to the existence of large solutions to the stationary equation
\begin{equation}\label{A3}
-\Gd w+f(w)=0,
\end{equation}
in $\Gw$, i.e. solutions which satisfy
\begin{equation}\label{A4}
\lim_{\gr_{_{\partial\Gw}}(x)\to 0}w(x)=\infty,
\end{equation}
and solutions of the ODE
\begin{equation}\label{A5}
\phi'+f(\phi)=0\quad\text{in }(0,\infty).
\end{equation}
subject to the initial blow-up condition
\begin{equation}\label{A6}
\lim_{t\to 0}\phi(t)=\infty.
\end{equation}
A natural assumption on $f$ is to assume that it is nondecreasing with $f(0)\geq 0$. If $f(a)>0$, a necessary and sufficient condition for the existence of a maximal solution $\overline w_{\Gw}$ to (\ref{A3}) is the Keller-Osserman condition, 
\begin{equation}\label{A7}
\myint{a}{\infty}\myfrac{ds}{\sqrt{F(s)}}<\infty,
\end{equation}
where $F(s)=\myint{0}{s}f(\gt)d\gt$. A necessary and sufficient condition for the existence of a solution $\phi$ of (\ref{A6}) with initial blow-up is 
\begin{equation}\label{A8}
\myint{a}{\infty}\myfrac{ds}{f(s)}<\infty.
\end{equation}
Furthermore the unique maximal solution $\overline\phi$ is obtained by inversion from the formula
\begin{equation}\label{A9}
\myint{\overline\phi(t)}{\infty}\myfrac{ds}{f(s)}=t\quad\forall t>0.
\end{equation}
It is known that, if $f$ is convex, (\ref{A7}) implies (\ref{A8}). If (\ref{A7}) holds and there exists a maximal solution to (\ref{A3}), it is not always true that this maximal solution is a large solution. In the case of a general nonlinearity, only sufficient conditions are known, independent of the regularity of $\prt\Gw$. We recall some of them. \smallskip

\noindent{\it If $N\geq 3$ and $f$ satisfies the weak singularity assumption
\begin{equation}\label{A10}
\myint{a}{\infty}s^{-2(N-1)/(N-2)}f(s)ds<\infty\quad\forall a>0.
\end{equation}
} \smallskip

\noindent{\it If $N=2$ and the exponential order of growth of $f$ defined by
\begin{equation}\label{A11}
a^+_{f}=\inf\left\{a\geq 0:\myint{0}{\infty}f(s)e^{-as}ds<\infty\right\}
\end{equation}
is finite.}

\medskip When $f(u)=u^q$ with $q>1$,  (\ref{A10}) means that $q<N/(N-2)$. When $q\geq N/(N-2)$ the regularity of $\prt\Gw$ plays a crucial role in the existence of large solutions. A necessary and sufficient condition involving a Wiener type test which uses the 
$C_{2,q'}^{\BBR^N}$-Bessel capacity has been obtained by probabilistic methods  by Dhersin and Le Gall \cite {DL} in the case $q=2$ and extended to the general case by Labutin \cite {La}. \medskip

 Uniqueness of the large solution of (\ref{A3}) has been obtained under three types of assumptions (see \cite {MV1}, \cite{MV4} and \cite{MV5}): \smallskip

\noindent{\it If $\prt\Gw=\prt\overline\Gw^c$ and $f(u)=u^q$ with $1<q<N/(N-2)$ or if $N=2$ and $f(u)=e^{au}$.}
\smallskip

\noindent{\it If $\prt\Gw$ is locally a continuous graph and $f(u)=u^q$ with $q>1$ or $f(u)=e^{au}$.}
 \smallskip

\noindent{\it If $f(u)=u^q$ with $q\geq N/(N-2)$ and $C^{\BBR^N}_{2,q'}
(\prt\Gw\setminus\tilde{\overline\Gw^c})=0$, where 
$\tilde E$ denotes the closure of a set in the fine topology associated to the Bessel capacity $C^{\BBR^N}_{2,q'}$.}
 \medskip
 
 In this article we extend most of the above mentioned results to the parabolic equation (\ref{A1}). We first prove that, {\it if $f$ 
 is super-additive, i. e. 
 \begin{equation}\label{A10'}
 f(x+y)\geq f(x)+f(y)\quad\forall (x,y)\in\BBR\ti\BBR,
\end{equation}
and satisfies (\ref{A7}) and (\ref{A8}), there exists a  maximal solution $\overline u_{Q^\Gw}$ to (\ref{A1}) in $Q^\Gw$, and it  satisfies
 \begin{equation}\label{A11'}
\overline u_{Q^\Gw}(x,t)\leq \overline w_{\Gw}(x)+\overline\phi(t)\quad\forall (x,t)\in Q^\Gw.
\end{equation}
If we assume also that $\prt\Gw=\prt\overline\Gw^c$, there holds}
 
  \begin{equation}\label{A12}
\max\{\overline w_{\Gw}(x),\overline\phi(t)\}\leq \overline u_{Q^\Gw}(x,t)\quad\forall (x,t)\in Q^\Gw.
\end{equation}

Under the assumption $\prt\Gw=\prt\overline\Gw^c$, it is possible to consider a decreasing sequence of smooth bounded domains $\Gw^{n}$ such that
$\overline\Gw^{n}\subset \Gw^{n-1}$, $\bar\Gw=\cap\Gw_n$, and prove that 
the increasing sequence of large solutions $\overline u_{Q^{\Gw^n}}$ of (\ref{A1}) in $Q^{\Gw^{n}}:=\Gw^{n}\ti (0,\infty)$,  converges to the {\it exterior maximal solution} $\underline u_{Q^\Gw}$ of (\ref{A1}) in $Q^\Gw$. If we proceed similarly with the large solutions $\overline w_{\Gw^{n}}$ of (\ref{A3}) in $\Gw^{n}$ and denote by 
$\underline w_{\Gw}$ their limit, then {\it  we prove that
  \begin{equation}\label{A13}
\max\{\underline w_{\Gw}(x),\overline\phi(t)\}\leq \underline u_{Q^\Gw}(x,t)\quad\forall (x,t)\in Q^\Gw.
\end{equation}}
The main result of this article is the following\medskip

\noindent{\bf Theorem 1. }{\it Assume $\Gw$ is a bounded domain such that $\prt\Gw=\prt\overline\Gw^c$, $f\in C(\BBR)$  is nondecreasing and satisfies (\ref{A7}), (\ref{A8}) and (\ref{A10'}). Then, if $\underline w_{\Gw}=\overline w_{\Gw}$, there holds $\underline u_{Q^\Gw}=\overline u_{Q^\Gw}$.}\medskip

Consequently, if (\ref{A3}) admits a unique large solution in $\Gw$, the same holds for (\ref{A1}) in $Q^\Gw_{\infty}$.
\mysection{The maximal solution}
In this section $\Gw$ is a bounded domain in $\BBR^N$ and $f\in C(\BBR)$ is nondecreasing and satisfies (\ref{A7}) and (\ref{A8}). We set $k_{0}=\inf\{\ell\geq 0:f(\ell)>0\}$ and assume also that, for any $m\in\BBR$ there exists $L=L(m)\in \BBR_{+}$ such that
  \begin{equation}\label{B1}
\forall (x,y)\in \BBR^2, x\geq m,\;y\geq m\Longrightarrow
f(x+y)\geq f(x)+f(y)-L.
\end{equation}
\bth{max} Under the previous assumptions there exists a maximal solution $\overline u_{Q^\Gw}$ in $Q^\Gw_{\infty}$.
\es
\Proof {\it Step 1- Approximation and estimates. } Let $\Gw_{n}$ be an increasing sequence of smooth domains such that $\overline\Gw_{n}\subset \Gw_{n+1}$ and $\cup\Gw_{n}=\Gw$. For each of these domains and $(n,k)\in \BBN^2_{*}$ we denote by $w=w_{n,k}$ the solutions of 
  \begin{equation}\label{B2}\left\{\BA {l}
-\Gd w+f(w)=0\quad\text{in }\Gw_{n}\\
\phantom{-\Gd +f(w)}
w=k\quad\text{in }\prt\Gw_{n}.
\EA\right.\end{equation}
where $\prt_{p}Q_{\infty}^{\Gw_{n}}:=\prt\Gw_{n}\ti (0,\infty)\cup 
\overline\Gw_{n}\ti\{0\}$. By \cite{Ke} there exists a decreasing function $g$ from $\BBR_{+}$ to $\BBR$, with limit $\infty$ at zero, such that 
  \begin{equation}\label{B4}
  w_{n,k}(x)\leq g\left(\gr_{_{\prt\Gw_{n}}}(x)\right)\quad\forall x\in\Gw_{n}.
\end{equation}
The mapping $k\to w_{n,k}$ is increasing, while $n\to w_{n,k}$ is decreasing. If we set
  \begin{equation}\label{B5}
\overline w_{\Gw}=\lim_{n\to\infty}\lim_{k\to\infty}  w_{n,k},
\end{equation}
it is classical that $\overline w_{\Gw}$ is the maximal solution of (\ref{A3}) in $\Gw$, and it satisfies
  \begin{equation}\label{B6}
  w(x)\leq g\left(\gr_{_{\prt\Gw}}(x)\right)\quad\forall x\in\Gw.
\end{equation}
We denote also by $u=u_{n,k}$ the solution of 
  \begin{equation}\label{B7}\left\{\BA {l}
u_{t}-\Gd u+f(u)=0\quad\text{in }Q_{\infty}^{\Gw_{n}}\\
\phantom{u_{t}-\Gd +f(u)}
u=k\quad\text{in }\prt_{p}Q_{\infty}^{\Gw_{n}}.
\EA\right.\end{equation}
By the maximum principle $k\to u_{n,k}$ is increasing and $n\to u_{n,k}$ decreasing. If we denote by $\bar\phi$ the maximal solution of the ODE (\ref{A5}), then $\bar\phi (t)$ is expressed by inversion by (\ref{A9}).
If 
$t_{k}= \bar\phi^{-1}(k)$, there holds, since $\bar\phi$ is decreasing, 
  \begin{equation}\label{B8}
\bar\phi(t+t_{k})\leq u_{n,k}(x,t)\quad\text {in }Q_{\infty}^{\Gw_{n}}.
\end{equation}
Furthermore, if $f(k)\geq 0$ (which holds if $k\geq k_{0}$), $w_{n,k}\leq k$. Therefore
  \begin{equation}\label{B9}
w_{n,k}(x)\leq u_{n,k}(x,t)\quad\text {in }Q_{\infty}^{\Gw_{n}}.
\end{equation}
Combining (\ref{B8}) and (\ref{B9}), we derive
  \begin{equation}\label{B10}
 \max\{w_{n,k}(x),\bar\phi(t+t_{k})\}\leq u_{n,k}(x,t)\quad\forall (x,t)\in Q_{\infty}^{\Gw_{n}}.
\end{equation}

Next we obtain an upper estimate. Let  $T>0$ and $m\in\BBR$ such that 
$$\min\{\overline w_{\Gw}(x):x\in\Gw\}> m\geq \bar\phi (T).$$ 
For $n\geq n_{1}$ and $k\geq k_{1}$ there holds
 $\min\{w_{n,k}(x):x\in\Gw\}\geq m$. Let $L=L(m)\geq 0$ be the corresponding damping term from (\ref{B1}). If $v_{n,k}=w_{n,k}(x)+\bar\phi(t+t_{k})$, then it satisfies
   \begin{equation}\label{B11}
v_{t}-\Gd v+f(v)=f(v)-f(\bar\phi(.+t_{k}))-f(w_{n,k})\geq -L
\quad\text{if }(x,t)\in \Gw_{n}\ti [0,T-t_{k}].
\end{equation}
Since $L\geq 0$, the function $\tilde v_{n,k}:=v_{n,k}+Lt$ is a supersolution for (\ref{A1}) in  $Q_{T-t_{k}}^{\Gw_{n}}:=\Gw_{n}\ti (0,T-t_{k})$ which dominates 
$u_{n,k}$ on $\prt_{p}Q_{T-t_{k}}^{\Gw_{n}}$, thus in $Q_{T-t_{k}}^{\Gw_{n}}$ by the maximum principle. Therefore
  \begin{equation}\label{B12}
 u_{n,k}(x,t)\leq w_{n,k}(x)+\bar\phi(t+t_{k})+Lt\quad\forall (x,t)\in Q_{T-t_{k}}^{\Gw_{n}}.
\end{equation}
{\it Step 2- Final estimates and maximality. } Using the different monotonicity properties of the mapping $(k,n)\mapsto w_{n,k}$ and the estimates (\ref{B10}) and (\ref{B12}), it follows that the function defined by
  \begin{equation}\label{B13}
\overline u_{Q^\Gw}:=\lim_{n\to\infty}\lim_{k\to \infty}u_{n,k}
\end{equation}
is a solution of (\ref{A1}) in $Q_{\infty}^{\Gw}$. Furthermore
  \begin{equation} \label{B14}
\max\{\overline w_{\Gw}(x),\bar\phi(t)\}\leq \overline u_{Q^{\Gw}}(x,t)\quad\forall (x,t)\in Q_{\infty}^{\Gw},
\end{equation}
and
  \begin{equation} \label{B15}
\overline u_{Q^{\Gw}}(x,t)\leq \overline w_{\Gw}(x)+\bar\phi(t)+tL(\phi (T)) \quad\forall (x,t)\in Q_{T}^{\Gw}.
\end{equation}
 since $\phi (T)\leq  \min\{\overline w_{\Gw}(x):x\in\Gw\}$. Next, we consider $u\in C^{2,1}(Q_{\infty}^\Gw)$, solution of (\ref{A1}) in $Q_{\infty}^\Gw$. Then, for $\ge>0$ and $n\in\BBN$, there exists $k^*>0$ such that for $k\geq k^*$, 
 $$ u_{n,k}(x,t-\ge)\geq u(x,t)\quad \forall (x,t)\in\Gw_{n}\ti (\ge,\infty).
 $$
 Letting successively $k\to\infty$, $n\to\infty$ and $\ge\to 0$, yields to $\overline u_{Q^\Gw}\geq u$ in $Q_{\infty}^\Gw$.\qeda\medskip
 
 Since $\overline w_{\Gw}$ be a large solution in $\Gw$ implies the same boundary blow-up for $\overline u_{Q^\Gw}$ on $\prt\Gw\ti (0,\infty)$, we give below some conditions which implies that $\overline u_{Q^\Gw}$ is a large solution.
 \bcor{large} Assume the assumptions of \rth{max} are fulfilled. Then $\overline u_{Q^\Gw}$ is a large solution if one of the following additional conditions is satisfied:\smallskip
 
 \noindent (i) $N\geq 3$ and $f$ satisfies the weak singularity condition (\ref{A10}).\smallskip
 
 \noindent (ii) $N=2$ and the exponential order of growth of $f$ defined by (\ref{A11}) is positive.\smallskip
 
 \noindent (iii) $N\geq 3$ and $\prt\Gw$ satisfies the Wiener regularity criterion. 
 \es
 \Proof Under condition (i) or (ii), for any $x_{0}\in \prt\Gw$, there exists a solution $w_{c,x_{0}}$ of 
   \begin{equation} \label{B16}\left\{\BA {l}
   -\Gd w+f(w)=c\gd_{x_{0}}\quad\text {in } B_{R}(x_{0})\\\phantom{   -\Gd +f(w)}
   w=0\quad\text {in } \prt B_{R}(x_{0}),
\EA\right.\end{equation}
where $R>0$ is chosen such that $\overline \Gw\subset B_{R}(x_{0})$ and $c>0$ is arbitrary under condition (i) and smaller that $2/a^+_{f}$ in case (ii). The function $w_{c,x_{0}}$ is radial with respect to $x_{0}$ and
$$\lim_{x\to x_{0}}w_{c,x_{0}}(x)=\infty.
$$
If $x\in \Gw$, we denote by $x_{0}$ a projection of $x$ on $\prt\Gw$. Since 
$$w_{n}(x)\geq w_{c,x_{0}}(x)\Longrightarrow
\overline w_{\Gw}(x)\geq w_{c,x_{0}}(x),$$
we derive from (\ref{B14}), 
$$\lim_{\gr_{_{\prt\Gw}}(x)\to 0}\overline u_{Q^\Gw}(x,t)=\infty,
$$
uniformly with respect to $t>0$.
In case (iii) we see that, for any $k>0$
   \begin{equation} \label{B17}\overline w_{\Gw}(x)\geq w_{k,\infty}(x)\quad\forall x\in \Gw,
\end{equation}
where $w_{k,\infty}$ is the solution of (\ref{B2}), with $\Gw_{n}$ replaced by $\Gw$. This again implies (\ref{B14}).\qeda\medskip

Using estimate (\ref{B14}) leads to the asymptotic behavior of $\overline u_{Q^\Gw}(x,t)$ when $t\to\infty$.

 \bcor{large-asym} Assume the assumptions of \rth{max} are fulfilled. Then $\overline u_{Q^\Gw}(x,t)\to \overline w_{\Gw}(x)$ locally uniformly on $\Gw$ when $t\to\infty$.
 \es
 \Proof For any $k>k_{0}$ and $n\in\BBN_{*}$ and any $s>0$, there holds by the maximum principle,
$$u_{n,k}(x,s)\leq k=u_{n,k}(x,0)\quad\forall x\in\Gw_{n}.$$
Using the monotonicty of $f$, we derive $u_{n,k}(x,t+s)\leq u_{n,k}(x,t)$ for any $(x,t)\in Q_{\infty}^{\Gw_{n}}$.
Letting $k\to\infty$ and then $n\to\infty$ yields to
   \begin{equation} \label{B18}\overline u_{Q^\Gw}(x,t+s)\leq \overline u_{Q^\Gw}(x,t)\quad\forall (x,t)\in Q_{\infty}^{\Gw}.   \end{equation}
It follows that $\overline u_{Q^\Gw}(x,t)$ converges to some $W(x)$ as $t\to\infty$ and $\overline w_{\Gw}\leq W$ from (\ref{B14}). Using the parabolic equation regularity theory, we derive that the trajectory $\CT:=\bigcup_{t\geq 0}\{\overline u_{Q^\Gw}(.,t)\}$ is compact in the $C^1_{loc}(\Gw)$-topology. Therefore $W$ is a solution of (\ref{A3}) in $\Gw$. It coincides with $\overline w_{\Gw}$ because of the maximality.\qeda
\mysection{Large solutions}

In this section we construct a minimal-maximal solution of (\ref{A1}) which is the minimal large solution whenever it exists. If $\prt\Gw$ is regular enough, the construction of the minimal large solution is easy.
\bth{minlarge}Let $\Gw$ be a bounded domain in $\BBR^N$ the boundary of which satisfies the Wiener regularity condition. If $f\in C(\BBR)$ is nondecreasing and satisfies (\ref{A7}), (\ref{A8}) and (\ref{B1}), then there exists a minimal large solution $\underline u_{Q^\Gw}$ to (\ref{A1}) in $Q^\Gw_{\infty}$. Furthermore 
\begin{equation}\label{C1}
\max\{\underline w_{\Gw}(x),\overline\phi(t)\}
\leq \underline u_{Q^\Gw}(x,t)
\quad\forall (x,t)\in Q^\Gw_{\infty},\end{equation}
and, for any $T>0$, 
\begin{equation}\label{C1'}
\underline u_{Q^\Gw}(x,t)\leq \underline w_{\Gw}(x)+\overline\phi(t)
+tL(\overline\phi (T))
\quad\forall (x,t)\in Q_{T}^\Gw,\end{equation}
where $L(\overline\phi (T))$ is as in (\ref{B17}), and  $\underline w_{\Gw}$ denotes the minimal large solution of (\ref{A3}) in $\Gw$.
\es
\Proof For $k\geq k_{0}$ (see Section 2), we denote by $\underline u_{k}$ the solution of  
\begin{equation}\label{C2}\left\{\BA {l}
u_{t}-\Gd u+f(u)=0\quad\text{in }Q^\Gw_{\infty}\\
\phantom{u_{t}-\Gd +f(u)}
u=k\quad\text{in }\prt_{p}Q^\Gw_{\infty}.
\EA\right.\end{equation}
When $k$ increases, $u_{k}$ increases and converges to some large solution $\underline u_{Q^\Gw}$ of (\ref{A1}) in $Q^\Gw_{\infty}$. If $u$ is any large solution of (\ref{A1}) in $Q^\Gw_{\infty}$, then the maximum principle and (\ref{A2}) implies $u\geq u_{k}$. Therefore 
$u\geq \underline u_{Q^\Gw}$. The same assumption allows to construct the solution $w_{k}$ of 
\begin{equation}\label{C3}\left\{\BA {l}
-\Gd w+f(w)=0\quad\text{in }\Gw\\
\phantom{-\Gd +f(w)}
w=k\quad\text{in }\prt\Gw,
\EA\right.\end{equation}
and, by letting $k\to\infty$, to obtain the minimal large solution $\underline w_{\Gw}$ of (\ref{A3}) in $\Gw$. Next we first observe, that, as in the proof of \rth{max}, (\ref{B11}) applies under the form
\begin{equation}\label{C4}
\overline\phi (t+t_{k})\leq u_{k}(x,t)\quad \text{in }Q^\Gw_{\infty},
\end{equation}
where, we recall it, $t_{k}=\overline\phi^{-1}(k)$. In the same way, for 
$k\geq k_{0}$ (with $f(k)\geq 0$), (\ref{B12}) holds under the form 
\begin{equation}\label{C5}
w_{k}(x)\leq u_{k}(x,t)\quad \text{in }Q^\Gw_{\infty}.
\end{equation}
Letting $k\to\infty$ yields to
\begin{equation}\label{C6}
\max\{\underline w_{\Gw}(x),\overline\phi(t)\}
\leq \underline u_{Q^\Gw}(x,t)
\quad\forall (x,t)\in Q^\Gw_{\infty}.\end{equation}

In order to prove the upper estimate we consider the same $m$ as it the proof of \rth{max} such that $\min\{\min\{w_{k}(x):x\in\Gw\},\overline\phi(t)\}\geq m$, and for $k'>k$, there holds
$$w_{k'}+\overline\phi \geq k=w_{k}\vline_{\prt_{p}Q^\Gw_{T}}.$$
Since $w_{k'}(x)+\overline\phi (t)+tL$ is a supersolution for (\ref{A1}) in $Q^\Gw_{T}$ it follows $w_{k'}+\overline\phi+tL \geq w_{k}$ in 
$Q^\Gw_{T}$. Letting successively $k'\to \infty$ and $k'\to \infty$, we derive (\ref{C1'}).\qeda

From this result we can deduce uniqueness results for solution of 

\bcor{uniq}Under the assumptions of \rth{minlarge}, if we assume moreover that $f$ is convex and, for any $\gth\in (0,1)$, there exists $r_{\gth}$ such that 
\begin{equation}\label{C6'}
r\geq r_{\gth}\Longrightarrow f(\gth r)\leq \gth f(r).
\end{equation}
Then
\begin{equation}\label{C7}
\underline w_{\Gw}=\overline w_{\Gw}\Longrightarrow
\underline u_{Q^\Gw}=\overline u_{Q^\Gw}.\end{equation}
\es
\Proof We fix $T\in (0,1]$ such that 
$$tL(\overline\phi (1))\leq \overline\phi (t)\quad\forall t\in (0,T],$$
 (remember that $L$ is always positive) and 
 $$2\underline w_{\Gw}(x)+\overline\phi(t)\geq 0\quad\forall (x,t)\in Q_{T}^\Gw.$$ 
 Then $\underline w_{\Gw}(x)+\overline\phi(t)\geq 0$ and
 $$\underline w_{\Gw}(x)+\overline\phi(t)+tL(\overline\phi (1))
 \leq\underline w_{\Gw}(x)+2\overline\phi(t) \leq\underline w_{\Gw}(x)+2\overline\phi(t) \leq 3\left(\underline w_{\Gw}(x)+\overline\phi(t)\right), $$
from which inequality follows
$$2^{-1}\left(\underline w_{\Gw}(x)+\overline\phi(t)\right)
\leq \underline u_{Q^\Gw}(x,t)\leq 3\left(\underline w_{\Gw}(x)+\overline\phi(t)\right)
\quad\forall (x,t)\in Q_{T}^\Gw.
$$
Therefore, if $\underline w_{\Gw}=\overline w_{\Gw}$, it follows
\begin{equation}\label{C8}
\underline u_{Q^\Gw}\leq \overline u_{Q^\Gw}\leq 6\underline u_{Q^\Gw}\quad\text{in } Q_{T}^\Gw.
\end{equation}
Next we assume $\underline u_{Q^\Gw}< \overline u_{Q^\Gw}$ and set
$$u^*=\underline u_{Q^\Gw}-
\myfrac{1}{6}\left(\overline u_{Q^\Gw}-\underline u_{Q^\Gw}\right).
$$
Since $f$ is convex, $u^*$ is a supersolution of (\ref{A1}) in $Q_{T}^\Gw$ (see \cite{MV2}, \cite{MV4}) and $u^*<\underline u_{Q^\Gw}$. Up to take a smaller $T$, we can also assume from (\ref{C6'}) that
$\min\{\underline u_{Q^\Gw}(x,t):(x,t)\in Q_{T}^\Gw\}\geq r_{1/12}$, thus
$$f(\underline u_{Q^\Gw}/12)\leq \myfrac{1}{12}f(\underline u_{Q^\Gw})\quad\text{in }Q_{T}^\Gw.
$$
Therefore $\underline u_{Q^\Gw}/12$ is a subsolution for (\ref{A1}) in $Q_{T}^\Gw$ and $12^{-1}\underline u_{Q^\Gw}< u^*$. Using a standard result of sub and super solutions and the fact that $f$ is locally Lipschitz continuous, we see that there exists some $u^{\#}$ solution of (\ref{A1}) in $Q_{T}^\Gw$ such that 
\begin{equation}\label{C9}
\myfrac{1}{12}\underline u_{Q^\Gw}\leq u^{\#}\leq u^*<\underline u_{Q^\Gw}\quad\text{in } Q_{T}^\Gw.
\end{equation}
Then $u^{\#}$ is a large solution, which contradicts the minimality of 
$\underline u_{Q^\Gw}$ on $Q_{T}^\Gw$. Finally 
$\underline u_{Q^\Gw}=\overline u_{Q^\Gw}$ in $Q_{\infty}^\Gw$.\qeda
\blemma {topo}Let $\Gw$ be a bounded domain in $\BBR^N$ and, for $\ge>0$, $\Gw_{\ge}:=\{x\in\BBR^N:\dist (x, \overline\Gw)<\ge\}$. The four following assertions are equivalent:\smallskip

\noindent (i) $\prt\Gw=\prt\overline\Gw^c$.\smallskip

\noindent (ii) For any $x\in\prt\Gw$, there exists a sequence $\{x_{n}\}\subset \overline\Gw^c$ such that $x_{n}\to x$.\smallskip

\noindent (iii) For any $x\in\prt\Gw$ and any $\ge>0$, 
$B_{\ge}(x)\cap \overline\Gw^c\neq\emptyset $.\smallskip

\noindent (iv) For any $x\in\prt\Gw$, $\lim_{\ge\to 0}\dist (x, \Gw_{\ge}^c)=0$.

\noindent (v) $\Gw=\overset{o}{\overline\Gw}$.
\es
\Proof There always holds $\prt\overline\Gw^c=\overline{\overline\Gw^c}\cap \overline\Gw
\subset \Gw^c\cap \overline\Gw=\prt\Gw$.\smallskip

\noindent (i)$\Longrightarrow$ (iii). Assume (iii) does not hold, there exist $x_{0}\in\prt\Gw$ and $\ge_{0}>0$ such that $B_{\ge_{0}}(x_{0})\cap\overline\Gw^c=\emptyset$. Thus $x_{0}\notin\overline{\overline\Gw^c}$, and $x_{0}\notin\prt\overline\Gw^c$. Therfore (i) does not hold.\smallskip

\noindent (iii)$\Longrightarrow$ (i). Let $x_{0}\in\prt\Gw$. If, for any 
$\ge>0$, $B_{\ge}(x)\cap\overline\Gw^c\neq\emptyset$, then
$x\in \overline{\overline\Gw^c}$. Because $x\in \Gw^c\cap \overline\Gw$, it implies that $x\in \overline\Gw\cap \overline{\overline\Gw^c}=\prt\overline\Gw^c$.
\smallskip

\noindent The equivalence between (iii) and (ii) is obvious.\smallskip

\noindent (ii))$\Longrightarrow$ (iv). We assume (iv) does not hold. There exist $x_{0}\in\prt\Gw$, $\ga>0$ and a sequence of positive real numbers $\{\ge_{n}\}$ converging to $0$ such that $\dist (x_{0},\Gw_{\ge_{n}}^c)\geq \ga$. Since for $\ge\geq\ge_{n}$, $\Gw_{\ge}^c\subset \Gw_{\ge_{n}}^c$, there holds $\dist (x_{0},\Gw_{\ge}^c)\geq \ga$. Furthermore, this inequality holds for any $\ge>0$. If there exist a sequence $\{x_{n}\}\subset \overline\Gw^c$ such that $x_{n}\to x_{0}$, then $\dist (x_{n},\overline\Gw)=\gd_{n}>0$, thus $x_{n}\in \Gw_{\gd_{n}}^c$. Consequently $|x_{n}-x_{0}|\geq\ga$, which is impossible. Therefore (ii) does not hold.\smallskip

\noindent (iv)$\Longrightarrow$ (iii). Let $x\in\prt\Gw$ and $x_{n}\in \Gw_{1\!/\!n}^c$ such that $|x-x_{n}|=\dist (x,\Gw_{1\!/\!n}^c)\to 0$. Since $\Gw_{1\!/\!n}^c\subset \overline\Gw$, $x_{n}\in \overline\Gw^c$ and $x_{n}\to x$.\smallskip

\noindent (iii)$\Longrightarrow$ (v). We first notice that $\overline\Gw=\cap_{\ge>0}\Gw_{\ge}=\cap_{\ge>0}\overline\Gw_{\ge}$ and
$\Gw\subset \overset{o}{\overline\Gw}$. If there exists some $x\in \overset{o}{\overline\Gw}\setminus\Gw$, then for some $\ge>0$, $B_{\ge}(x)\subset\overline\Gw$ which implies $B_{\ge}(x)\cap \overline\Gw^c=\emptyset$. But $x\notin\Gw$ implies $x\in\prt\Gw$. Thus (iii) does not hold.\smallskip

\noindent (v)$\Longrightarrow$ (iii).  If (iii) does not hold, there exists 
$x\in \prt\Gw$ and $\ge>0$ such that $B_{\ge}(x)\cap\overline \Gw^c=\emptyset\Longleftrightarrow B_{\ge}(x)\subset \overline \Gw$. Therefore $x\in \overset{o}{\overline\Gw}\setminus\Gw$.
\qeda\medskip
\bdef{min-max}A solution $U$ (resp. $W$ to problem (\ref{A1}) in  $Q_{\infty}^\Gw$ (resp. (\ref{A3}) in $\Gw$) is called an exterior maximal solution if it is larger than the restriction to $Q_{\infty}^\Gw$
(resp. $\Gw$) of any solution of (\ref{A1}) (resp. (\ref{A3}) ) defined in an open neighborhood of $Q_{\infty}^\Gw$ (resp.  $\Gw$)).\es
\bprop{minmax-E} Assume $\Gw$ is a bounded domain in $\BBR^N$ such that $\prt\Gw=\prt\overline\Gw^c$ and  $f\in C(\BBR)$ is nondecreasing and satisfies (\ref{A7}). Then there exists an exterior maximal solution $\underline w_{\Gw}^*$ to problem (\ref{A3}) in $\Gw$.
\es
\Proof Since $\prt\Gw=\prt\overline\Gw^c$ we can consider the decreasing sequence of  the $\Gw_{1\!/\!n}$ defined in \rlemma {topo} with $\ge=1\!/\!n$ and, for each $n$, the minimal large solutions $\underline w_{n}$ of (\ref{A3}) in $\Gw_{1\!/\!n}$: this possible since $\prt\Gw_{1\!/\!n}$ is Lipschitz. The sequence $\{\underline w_{n}\}$ is increasing. Its restriction to $\Gw$ is bounded from above by the maximal solution $\overline w_{\Gw}$. It converges to some function 
$\underline w_{\Gw}^*$. By \rlemma{topo}-(v), $\underline w_{\Gw}^*$ is a solution of (\ref{A3}) in the interior of $\cap_{n}\Gw_{1\!/\!n}$ which is $\Gw$. If $w$ is any solution of (\ref{A3}) defined in an open neighborhood of $\overline\Gw$, it is defined in $\Gw_{1\!/\!n}$ for $n$ large enough and therefore smaller than $\underline w_{n}$. Thus $w\vline_{\Gw}\leq \underline w_{\Gw}^*$. 
Consequently, $ \underline w_{\Gw}^*$ coincides with the supremum of the restrictions to 
$\Gw$ of solutions of (\ref{A3}) defined in an open neighborhood of $\overline \Gw$.
\qeda
\bprop {W} Let $f\in C(\BBR)$ be a nondecreasing function for which (\ref{A7}) holds and $\Gw$ a bounded domain in $\BBR^N$ such that $\prt\Gw=\prt\overline\Gw^c$. Then $ \underline w_{\Gw}^*$ is smaller than any large solution. Furthermore, if $\prt\Gw$ satisfies the Wiener regularity criterion and is locally the graph of a continuous function, then $\underline w_{\Gw}=\underline w_{\Gw}^*$. 
\es
\Proof  We first notice that  Wiener criterion implies  statement (iii) in \rlemma{topo}, hence $\prt\Gw=\prt\overline\Gw^c$. If $w_{\Gw}$ is a large solution, it dominates on $\prt\Gw$, and therefore in $\Gw$ by the maximum principle, the restriction to $\Gw$ of any function $w$ solution of (\ref{A3}) in an open neighborhood of $\overline\Gw$. Then 
$$\underline w_{\Gw}^*\leq  w_{\Gw}.$$
Consequently, if $\underline w_{\Gw}^*$ is a large solution, it coincides with the minimal large solution $\underline w_{\Gw}$. Because $\prt\Gw$ is compact, there exists a finite number of bounded open subset $\CO_{j}$, hyperplanes $H_{j}$ and continuous functions $h_{j}$ from $H_{j}\cap \overline \CO_{j}$ into $\BBR_{+}$ such that 
$$\prt\Gw\cap\overline \CO_{j}=\left\{x=x'+h_{j}(x')\gn_{j}:\;\forall x'\in H_{j}\cap \overline \CO_{j}\right\}
$$
where $\gn_{j}$ is a fixed unit vector orthogonal to $H_{j}$ and $\prt\Gw\subset\cup_{j}\CO_{j}$. We can assume that $H_{j}\cap \overline \CO_{j}=\overline B_{j}$ is a (N-1) dimensional closed ball and, 
$$G_{j}:=\{x=x'+t\gn_{j}:\;x'\in \overline B_{j},\,0\leq t< h_{j}(x')\}\subset \Gw,$$
$$G^{\#}_{j}:=\{x=x'+t\gn_{j}:\;x'\in \overline B_{j},\,h_{j}(x')< t\leq a \}\subset \overline\Gw^c.,$$
for some $a>0$ such that $a/4<h_{j}(x')<3a/4$ for any $x'\in \overline B_{j}$. Finally, we can assume that 
$$\CO_{j}=\{x=x'+t\gn_{j}:\;x'\in \overline B_{j},\, 0\leq t\leq a\}.$$
Let $\ge\in (0,a/8)$ and
$$G_{j,\ge}:=\{x=x'+t\gn_{j}:\;x'\in \overline B_{j},\,\ge\leq t< h_{j}(x')+\ge\}.
$$
There exists a smooth bounded domain $\Gw'$ such that $\overline\Gw\subset \Gw'$ and
$$\prt\Gw'\cap \overline\CO_{j}=\{x=x'+\ell (x')\gn_{j}:x'\in\overline B_{j},\,
h(x')+\ge/2\leq \ell(x')\leq h(x')+3\ge/2\},$$
where $\ell\in C^{\infty}(\overline B_{j})$. We denote $G_{j}:=G_{j,0}$,
$$\prt_{p}G_{j,\ge}:=\{x=x'+t\gn_{j}:\;x'\in \prt B_{j},\,\ge\leq t\leq h_{j}(x')+\ge\}\cup \{x=x'+ \ge\gn_{j}:\;x'\in  B_{j} \},$$ 
and
$$\prt_{u}G_{j,\ge}:=\{x=x'+(h_{j}(x')+\ge)\gn_{j}:\;x'\in  B_{j} \}.
$$
Let $w'$ be the minimal large solution of (\ref{A3}) in $\Gw'$, $\ga'=\min \{w'(x):x\in \Gw'\}$ and  $W_{\ge}$  the minimal solution of 
\begin{equation}\label{}\left\{\BA {l}
-\Gd  W+f(W)=0\quad\text{in }G_{j,\ge}\\
\phantom{-\Gd  +f(W)}
W=\ga'\quad\text{in }\prt_{p}G_{j,\ge}\\
\displaystyle\lim _{t\to h(x')+\ge}W(x'+t\gn_{j})=\infty\quad\forall x'\in  B_{j}.
\EA\right.\end{equation}
Then $w'\geq W_{\ge}$ in $G_{j,\ge}\cap \Gw'$. Furthermore
$W_{\ge}(x)=W_{\ge}(x'+t\gn_{j})=W_{0}(x'+(t-\ge)\gn_{j})$ for any $x'\in \overline B_{j}$ and $\ge<t<h(x')+\ge$. Therefore, given $k>0$, there exists $\gd_{k}>0$ such that for any 
$$x'\in \overline B_{j}\text { and }h_{j}(x')-\gd_{k}\leq t<h_{j}(x')\Longrightarrow W_{0}(x'+t\gn_{j})\geq k.
$$
As a consequence, 
$\liminf_{t\to h_{j}(x')} \underline w_{\Gw}^*(x'+t\gn_{j})\geq k$, uniformly  with respect to $x'\in\overline B_{j}$. This implies that $\underline w_{\Gw}^*$ is a large solution.\qeda\medskip

\noindent\Remark We conjecture that the equality  $\underline w_{\Gw}^*=\underline w_{\Gw}$ holds under the mere assumption that the Wiener criterion is satisfied. 

\bth{min-maxth} Assume $\Gw$ is a bounded domain in $\BBR^N$ such that $\prt\Gw=\prt\overline\Gw^c$ and $f\in C(\BBR)$ satisfies (\ref{A7}), (\ref{A8}) and (\ref{B1}). Then there exists a exterior maximal solution $\underline u_{Q^\Gw}^*$ to problem (\ref {A1}). Furthermore estimates (\ref{C1}) and (\ref{C1'}) hold with $\underline w_{\Gw}$ replaced by the exterior maximal solution $\underline w_{\Gw}^*$ to problem (\ref{A3}) in $\Gw$.
\es
\Proof The construction of $\underline u_{Q^\Gw}^*$ is similar to the one of $\underline w_{\Gw}$, since we can restrict to consider open neighborhoods $Q_{1\!/\!n}= \Gw_{1\!/\!n}\ti (-1\!/\!n,\infty)$. Then 
$\underline u_{Q^\Gw}^*$ is the increasing limit of the minimal large solutions $u_{n}$ of (\ref{A1}) in $Q_{1\!/\!n}$, since $\overline{Q_{\infty}^{\Gw}}=\cap_{n}Q_{1\!/\!n}$ and, by \rlemma{topo}-(v), $Q_{\infty}^{\Gw}=\overset{o}{\overline {Q_{\infty}^{\Gw}}}$. We recall that the minimal large solution $w_{n}$ of (\ref{A3}) in $\Gw_{1\!/\!n}$ is the increasing limit, when $k\to\infty$, of the sequence of solution $\{w_{n}^k\}$ of 
\begin{equation}\label {C12'}\left\{\BA {l}
-\Gd w+f(w)=0\quad\text {in }\Gw_{1\!/\!n}\\[2mm]
\phantom{-\Gd w+f()}
w=k\quad\text {on }\prt\Gw_{1\!/\!n},
\EA\right.\end{equation}
while the minimal large solution $u_{n}$ of (\ref{A1}) in $Q_{1\!/\!n}$ is the (always increasing) limit of the solutions $u_{n}^k$ of 
\begin{equation}\label {C12}\left\{\BA {l}
u_{t}-\Gd u+f(u)=0\quad\text {in }Q_{1\!/\!n

}\\[2mm]
\phantom{u_{t}-\Gd u+f()}
u=k\quad\text {on }\prt_{p}Q_{1\!/\!n}.
\EA\right.\end{equation}
Clearly
$$ \max\{w_{n}^k,\overline\phi (.+1\!/\!n)\}\leq u_{n}(x,t),
$$
which implies (\ref{C1}). For the other inequality, we see that 
$(x,t)\mapsto w_{n}^k(x)+\overline\phi (t)+Lt$ is a supersolution which dominates $u_{n}^k$ on $\prt_{p}$, where $L$ corresponds to the minimum of $w_{n}^k$ in $\Gw_{1\!/\!n}Q_{1\!/\!n}$. Thus
$$u_{n}(x,t)\leq w_{n}^k+\overline\phi (.+1\!/\!n),
$$
which implies 
\begin{equation}\label {C13}
\max\{\underline w_{\Gw}^*(x),\overline\phi (t)\}\leq \underline u_{\Gw}^*(x,t)\quad\forall (x,t)\in Q_{\infty}^\Gw.
\end{equation}
The upper estimate is proved in the following way. If $k>n$, $\overline Q_{k}\subset Q_{n}$. Therefore, choosing $m$ such that 
$\min\left\{\min\{\underline w_{\Gw_{1\!/\!k}}(x):x\in \Gw_{1\!/\!k} ,\min\{\overline \phi(t+1\!/\!k):t\in (0,T]\}\right\}\geq m$, we obtain that 
$(x,t)\mapsto \underline w_{\Gw_{1\!/\!k}}(x)+\phi(t+1\!/\!k)+Lt$ is a super solution of (\ref{A1}) in $Q_{T}^{\Gw_{1\!/\!k}}$, thus it dominates the minimal large solution of (\ref{A1}) in $Q_{T}^{\Gw_{1\!/\!n}}$. Letting successively $k\to\infty$ and $n\to\infty$, yields to 
\begin{equation}\label {C14}
\underline u_{\Gw}^*(x,t)\leq \underline w_{\Gw}^*(x)+\overline\phi (t)\quad\forall (x,t)\in Q_{T}^\Gw.
\end{equation}
\qeda\noindent

The next result extends \rcor{uniq} without the boundary Wiener regularity assumption.

\bth{uniq2} Let $\Gw$ be a bounded domain in $\BBR^N$ such that $\prt\Gw=\prt\overline\Gw^c$. If $f\in C(\BBR)$ is convex and satisfies (\ref{A7}), (\ref{A8}), (\ref{B1}) and (\ref{C6'}). Then, if $\underline w_{\Gw}^*$ is a large solution, the following implication holds
\begin{equation}\label{C15}
\underline w_{\Gw}^*=\overline w_{\Gw}\Longrightarrow
\underline u_{Q^\Gw}^*=\overline u_{Q^\Gw}.\end{equation}
\es
\Proof If $\underline w_{\Gw}^*$ is a large solution, the same is true for $\underline u_{Q^\Gw}^*$ because of (\ref{C1}). Actually $\underline u_{Q^\Gw}^*$ is the minimal large solution in $Q^\Gw_{\infty}$ for the same reasons as $\underline w_{\Gw}^*$. Therefore the proof of \rcor{uniq} applies and it implies the result. \qeda\\

\noindent \Remark We conjecture that (\ref{C15}) holds, even if $\underline w_{\Gw}^*$ is not a large solution.

 

\begin{thebibliography}{99}
 \bibitem{BP}P. Baras \& M. Pierre, {\em Probl\`emes paraboliques semi-lin\'eaires avec
donn\'ees mesures}, Applicable Anal. {\bf 18}, 111-149 (1984).

\bibitem{BPT}  H. Brezis, L. A. Peletier \& D. Terman, {\em A very singular solution of the heat equation with absorption}, Arch. Rat. Mech. Anal. {\bf 95}, 185-209 (1986).

\bibitem{BF}  H. Brezis and A. Friedman, {\em Nonlinear parabolic equations involving
measures as initial conditions}, J. Math. Pures Appl. {\bf 62}, 73-97 (1983).


\bibitem{DL} J. S. Dhersin and J. F. Le Gall, \textit{ Wiener's test for super-Brownian motion and the Brownian snake}, Probab. Theory Relat. Fields {\bf 108}, 103-29 (1997).

\bibitem{Ke} J.B. Keller,\textit{ On solutions of $\Delta u=f(u)$}, Comm. Pure Appl. Math. {\bf 10}, 
503-510 (1957).

\bibitem{La} D. Labutin, \textit{ Wiener regularity for large solutions of nonlinear equations}, Archiv f\"or Math. {\bf 41}, 307-339 (2003).

\bibitem{MV1} M.  Marcus and L. V\'eron, \textit{ Uniqueness and
asymptotic behaviour of solutions with boundary blow-up for a class of nonlinear elliptic equations},  Ann. Inst. H. Poincar\'e {\bf 14}, 237-274 (1997).


\bibitem{MV2}  M.  Marcus and L. V\'eron, \textit{ The boundary trace of positive solutions of semilinear elliptic equations: the subcritical case}, Arch. Rat. Mech. Anal. {\bf 144}, 201-231 (1998).

\bibitem{MV3}  M.  Marcus and L. V\'eron, \textit{ The initial trace of positive solutions of semilinear parabolic equations}, Comm. Part. Diff. Equ. {\bf 24}, 1445-1499 (1999).

\bibitem{MV4}  M.  Marcus and L. V\'eron, \textit{ Existence and uniqueness results for large solutions of general
nonlinear elliptic equations}, J. Evol. Equ {\bf 3}, 637-652 (2003).

\bibitem{MV5}  M.  Marcus and L. V\'eron, \textit{Maximal solutions for $-\Delta u+u^q=0$ in open and finely open sets},   J. Math. Pures Appl. {\bf 91}, 256-295 (2009).

\bibitem{Oss} R. Osserman, \textit{ On the inequality $\Delta u\geq f(u)$}, Pacific J. Math.{\bf 7}, 1641-1647  (1957).

\bibitem{Ve1} L. V\'eron, \textit{ Generalized boundary value problems for nonlinear elliptic equations}, Electr. J. Diff. Equ. Conf. {\bf 6}, 313-342 (2000).

\bibitem{Ve2} L. V\'eron, \textit{ Singular solutions of some nonlinear elliptic equations}, Nonlinear Anal. T. M. \& A {\bf  5}, 225-242 (1981).

 \end{thebibliography}
 \end{document}